\documentclass[11pt,a4paper]{article}

\usepackage{amsmath}
\usepackage{amsthm}
\usepackage{amssymb}
\usepackage{amscd}

\title{On Pseudo-Effectivity of the Second Chern Classes for Terminal Threefolds}
\author{Qihong Xie}
\date{}
\theoremstyle{plain}
\newtheorem{prop}{Proposition}[section]
\newtheorem{lem}[prop]{Lemma}
\newtheorem{thm}[prop]{Theorem}
\newtheorem{cor}[prop]{Corollary}
\newtheorem{conj}[prop]{Conjecture}
\newtheorem{prob}[prop]{Problem}

\theoremstyle{definition}
\newtheorem{defn}[prop]{Definition}
\newtheorem*{ack}{Acknowledgment}

\theoremstyle{remark}
\newtheorem{rem}[prop]{Remark}
\newtheorem{ex}[prop]{Example}

\newcommand{\Q}{\mathbb Q}
\newcommand{\R}{\mathbb R}

\newcommand{\Z}{\mathbb Z}
\newcommand{\N}{\mathbb N}

\newcommand{\PP}{\mathbb P}
\newcommand{\OO}{\mathcal O}
\newcommand{\II}{\mathcal I}

\newcommand{\NN}{\mathcal N}

\newcommand{\TT}{\mathcal T}

\newcommand{\A}{\mathcal A}

\newcommand{\Pic}{\mathop{\rm Pic}\nolimits}
\newcommand{\Supp}{\mathop{\rm Supp}\nolimits}

\begin{document}

\maketitle

\begin{abstract}
We give a reduction of the conjecture that for terminal projective threefolds whose anticanonical divisors are nef, the second Chern classes are pseudo-effective. On the other hand, some effective non-vanishing results are obtained as applications of the pseudo-effectivity of the second Chern classes.
\end{abstract}

\setcounter{section}{0}
\section{Introduction}\label{C:S1}

The second Chern class plays an essential role in the three-dimensional birational geometry. For example, as a famous result, the Miyaoka theorem says that the second Chern class is pseudo-effective for a terminal projective minimal threefold (cf.\ \cite{mi}). 
Furthermore, by the pseudo-effectivity of the second Chern class, 
we can prove the non-negativity of the Kodaira dimension and 
the abundance theorem for a terminal projective minimal threefold.

In this paper, we consider the following:

\begin{conj}\label{C:1.1}
Let $X$ be a terminal projective threefold with $-K_X$ nef. Then the second Chern class $c_2(X)$ is pseudo-effective.
\end{conj}

The second Chern class of a terminal threefold $X$, 
as a 1-cycle, is defined as follows: 
$c_2(X):=c_2(\TT_X|_U)\in A^2(U)\cong A^2(X)$, 
where $U$ is the smooth locus of $X$ and $\TT_X$ is the tangent sheaf on $X$. 
By definition, a 1-cycle is said to be pseudo-effective, if its numerical equivalence class is contained in the Kleiman-Mori cone $\overline{NE}(X)$. For Conjecture \ref{C:1.1}, we have the following known results (cf.\ \cite{xie}, Theorem 2.2 and Proposition 2.4).

\begin{thm}\label{C:1.2}
Let $X$ be a terminal projective threefold with $-K_X$ nef. If the numerical dimension $\nu(-K_X)\neq 2$, then $c_2(X)$ is pseudo-effective. If $\nu(-K_X)=2$, then $c_1(X).c_2(X)\geq 0$ and the irregularity $q(X)\leq 1$.
\end{thm}

We say that an extremal contraction $f: X\rightarrow Y$ is good, if there exists an integer $n\geq 0$ such that $c_2(X)+nl$ is pseudo-effective, where $\R_+[l]$ is the corresponding extremal ray of $f$.

As for Conjecture \ref{C:1.1}, we have considered a simpler case when $X$ is smooth in \cite{xie}. As a result, a proof has been given for the smooth case under a weak assumption ${\rm (AD_{III})}$. The main idea of the proof is to investigate the goodness of all extremal contractions from $X$, which implies the pseudo-effectivity of $c_2(X)$. If $q(X)=1$, a nice classification is available, which guarantees a complete proof. The case $q(X)=0$ is more complicated, because not only a similar nice classification is unavailable, but also in the subcase ${\rm (D_{III})}$, we cannot prove that $f$ is good by induction. Therefore, ${\rm (AD_{III})}$ is such an assumption that the subcase ${\rm (D_{III})}$ is good.

In order to rule out the assumption ${\rm (AD_{III})}$, or say, to run the Minimal Model Program in the subcase ${\rm (D_{III})}$, it is necessary to consider the class of $\Q$-factorial terminal projective threefolds with almost nef anticanonical divisors, instead of the class of smooth projective threefolds with nef anticanonical divisors. In other words, it is better to start from the terminal case than the smooth case to prove Conjecture \ref{C:1.1}.

We use the same idea to extend the argument of smooth case to terminal case. It is obvious that the terminal case is more delicate than the smooth case, since we have to consider the flips when running the MMP. But it turns out that the flips behave well in this new class. Unfortunately, since in general, it is difficult to bound the terminal singularities when running the MMP, we cannot obtain a complete proof of the terminal case without an assumption on the non-negativity of $c_1.c_2$.

\begin{defn}\label{C:1.3}
Let $X$ be a $\Q$-factorial terminal projective threefold. A threefold $Y$ is in the class $\A(X)$, if there is a composition of birational maps: $X=X_0\dashrightarrow X_1\dashrightarrow\cdots\dashrightarrow X_r=Y$, where $\dashrightarrow$ is either a divisorial contraction or a flip.
\end{defn}

The following is the main theorem in this paper.

\begin{thm}\label{C:1.4}
Let $X$ be a terminal projective threefold such that $-K_X$ is nef and $\nu(-K_X)=2$.

{\rm (1)} If $q(X)=1$, then $c_2(X)$ is pseudo-effective;

{\rm (2)} If $q(X)=0$, then $c_2(X)$ is pseudo-effective provided that $c_1(Y).c_2(Y)\geq 0$ holds for any $Y\in\A(X^\Q)$, where $\mu: X^\Q\rightarrow X$ is a $\Q$-factorialization of $X$.
\end{thm}

In \S \ref{C:S2}, we give a complete proof of Conjecture \ref{C:1.1} when $q(X)=1$. In \S \ref{C:S3}, a reduction of Conjecture \ref{C:1.1} is proved when $q(X)=0$. In \S \ref{C:S4}, we give some applications to the Effective Non-vanishing Conjecture.

In the whole paper, we will use freely the results on the Minimal Model Theory from \cite{kmm} and \cite{km}. For some necessary definitions and notation, we refer to \cite{xie}.

We work over the field of complex numbers.

\begin{ack}
I am very grateful to Professor Yujiro Kawamata for his valuable advice and warm encouragement. I would also like to thank Professors Keiji Oguiso and Hiromichi Takagi, and Doctors Masayuki Kawakita, Yasunari Nagai and Shunsuke Takagi for stimulating discussions.
\end{ack}

\section{Proof of the case $q(X)=1$}\label{C:S2}

First, we give a lemma, whose proof is similar to that of Proposition 3.3 of \cite{dps}.

\begin{lem}\label{C:2.1}
Let $X$ be a Gorenstein $\Q$-factorial terminal projective threefold with $-K_X$ nef. Let $f: X\rightarrow Y$ be an extremal contraction which contracts a divisor $E$ to a curve $C$. If $-K_Y$ is not nef, then $C$ is a rational curve.
\end{lem}

\begin{proof}
It follows from Lemmas 2 and 3 of \cite{cu} that both $X$ and $Y$ are factorial, and $f$ is just the blow-up of $Y$ along $C$, which is a local complete intersection. From the formula $K_X=f^*K_Y+E$, we immediately see that $-K_Y.C'\geq 0$ for every curve $C'\neq C$ in $Y$.

Let $\NN_{C|Y}=\NN$ be the normal bundle of $C$ to $Y$. Since $\NN$ is locally free, $E\cong\PP(\NN^*)$ is a $\PP^1$-bundle over $C$. Let $\pi: E\rightarrow C$ be the projection onto $C$, $F$ the general fiber of $\pi$. Since $E$ is Cartier, $K_E=(K_X+E)|_E$. It follows from $K_E.F=-2$ that $K_X.F=E.F=-1$. Since $\Pic(E)\cong\pi^*\Pic(C)\oplus\Z$, there is a Cartier divisor $C_1$ on $E$ such that $\Pic(E)$ is numerically generated by $C_1$ and $F$, and $C_1.F=1$. It is easy to see that $C_1$ is irreducible and reduced, and $\pi|_{C_1}: C_1\rightarrow C$ is a birational surjective morphism. Assume that $e=-C_1^2$ and $p_a(C_1)$ is the arithmetical genus of $C_1$ defined by $2p_a(C_1)-2=C_1.(C_1+K_E)$. Since $K_X.F=-1$, we may assume that
\[ -K_X|_E\equiv C_1+bF, \]
where $b\in\Z$. Since $\OO_E(E)\cong \OO_E(-1)$, we may assume that
\[ E|_E\equiv -C_1+\mu F, \]
where $\mu\in\Z$. Hence $(-K_Y.C)=(-K_X.C_1)+(E.C_1)=b+\mu$, so $-K_Y$ is nef if $b+\mu\geq 0$. It is easy to see that $K_E\equiv -2C_1+(2p_a(C_1)-2-e)F$, hence $K_E^2=8(1-p_a(C_1))$. On the other hand, $K_E^2=(K_X|_E+E|_E)^2=(-2C_1+(\mu-b)F)^2$. Hence $b+\mu=2b-e+2(p_a(C_1)-1)$.

Since $-K_X|_E$ is nef, we have that $(-K_X|_E)^2\geq 0$, which implies that $b\geq e/2$. If $p_a(C_1)\geq 1$, then $-K_Y$ is nef. Thus $p_a(C_1)=0$, namely $C_1$ is a smooth rational curve. Hence $C$ is rational.
\end{proof}

\begin{prop}\label{C:2.2}
Let $X$ be a $\Q$-factorial terminal projective threefold such that $-K_X$ is nef, $\nu(-K_X)=2$ and $q(X)=1$. Let $f: X\rightarrow Y$ be an extremal contraction. Then $(Y,f)$ is one of the following cases.

${\rm (F_I)}$ $X$ is smooth, $Y$ is a smooth elliptic curve and $f$ is a del Pezzo fibration.

${\rm (F_{II})}$ Let $\alpha: X\rightarrow A$ be the Albanese map of $X$. Then $X$ is smooth, $Y$ is a smooth hyperelliptic surface or a $\PP^1$-bundle over $A$ with $-K_Y$ nef. $f$ is a conic bundle.

{\rm (D)} $f:X\rightarrow Y$ is a divisorial contraction which contracts a divisor $E$ to a curve $C$. Both $X$ and $Y$ are Gorenstein, $-K_Y$ is nef and $\nu(-K_Y)\leq 2$. $f$ is just the blow-up of $Y$ along $C$.
\end{prop}

\begin{proof}
It is easy to see that $\dim Y>0$. Since $-K_X$ is nef but not big, and $\chi(\OO_X)=1-q(X)=0$, we have that $X$ is Gorenstein (cf.\ \cite{matsuki}, Corollary 6.3), hence factorial.

(\ref{C:2.2}.1) If $\dim Y=1$, then $Y$ is a smooth elliptic curve since $q(Y)=q(X)=1$, and $f$ is a del Pezzo fibration. The smoothness of $X$ follows from Proposition 1.3 of \cite{ps}.

(\ref{C:2.2}.2) If $\dim Y=2$, then $Y$ is a smooth surface and $f$ is a conic bundle (cf.\ \cite{cu}, Theorem 7). The explicit structure of $Y$ is determined by Proposition 1.7 of \cite{ps}.

(\ref{C:2.2}.3) Assume that $\dim Y=3$. By the structure theorem for $f$ (cf.\ \cite{cu}), $f$ could never be a small contraction. If $f$ contracts a divisor $E$ to a terminal point $p\in Y$, then $-K_Y$ is also nef, furthermore big through a simple computation. It follows from the Kawamata-Viehweg vanishing theorem that $q(X)=q(Y)=0$, it is absurd. Hence $f$ contracts a divisor $E$ to a possibly singular curve $C$, which is a locally complete intersection in $Y$. Note that $f$ is just the blow-up of $Y$ along $C$ and that $Y$ is also factorial.

Assume that $-K_Y$ is not nef, then $C$ is a rational curve by Lemma \ref{C:2.1}. We may repeat the same argument as in the proof of Theorem 3.3(3) of \cite{ps} to obtain that $C$ is a multi-section of the Albanese map $\alpha: Y\rightarrow {\rm Alb}(Y)$, hence $C$ is a smooth elliptic curve. This is a contradiction. Thus $-K_Y$ is nef. It is easy to see that $\nu(-K_Y)\leq 2$.
\end{proof}

\begin{lem}\label{C:2.3}
Let $X$ be a terminal projective threefold, $\R_+[l]$ an extremal ray on X, and $\alpha$ a positive number. Assume that $c_2(X)+nl$ is pseudo-effective for some $n\in\N$. Then we can take an ample Cartier divisor $L$ on $X$, a sufficiently small $\varepsilon>0$, and the cone decomposition $\overline{NE}(X)=\R_+[l]+\sum_i\R_+[l_i]+\overline{NE}_\varepsilon(X)$ such that for any decomposition $c_2(X)=al+\sum_ib_il_i+z$, where $b_i\geq 0$, $z\in\overline{NE}_\varepsilon(X)$, we have $z.(-K_X)<\alpha$.
\end{lem}

\begin{proof}
The proof is similar to that of Proposition 4.10 of \cite{xie}.
\end{proof}

\begin{thm}\label{C:2.4}
Let $X$ be a terminal projective threefold such that $-K_X$ is nef, $\nu(-K_X)=2$ and $q(X)=1$. Then $c_2(X)$ is pseudo-effective.
\end{thm}

\begin{proof}
First, we assume that $X$ is $\Q$-factorial. In cases ${\rm (F_I)}$ and ${\rm (F_{II})}$, since $X$ is smooth, $c_2(X)$ is pseudo-effective by Theorem 3.12 of \cite{xie}. In case (D), we use induction on the Picard number $\rho(X)$ to prove the pseudo-effectivity of $c_2(X)$.

With notation and assumptions as in case (D). Assume that $c_2(Y)$ is pseudo-effective. Let $S=\{p\in Y|f^{-1}(p)$ contains some singular point of $X$\}. Then $S$ is a finite set of points and $f: X\setminus f^{-1}(S)\rightarrow Y\setminus S$ is a blow-up along a smooth curve $C$. Let $F$ be a general fiber of $f|_E: E\rightarrow C$. Then
\[ c_2(X)=f^*c_2(Y)+C_1+xF,\]
where $C_1$ is the curve in $E$ as in Lemma \ref{C:2.1}, and $x$ is a rational number (cf.\ \cite{xie}, Proposition 3.8). By the same reason, we have
\[ f^*C=C_1+yF \]
where $y$ is a rational number.

Let $n_1$ (resp.\ $n_2$) be 0 when $x$ (resp.\ $y$) is non-negative, otherwise the smallest integer not less than $-x$ (resp.\ $-y$). Assume that $c_2(Y)=\lim_{k\to\infty}\xi_k=\lim_{k\to\infty}(a_kC+R_k)$, where $\xi_k$ are effective and $C$ is not contained in the support of $R_k$. Then there is a positive integer $N$ such that $\sup_k\{a_k\}<N$. Let $H$ be a nef Cartier divisor on $X$. Then $f_*H$ is nef except along $C$.
\begin{eqnarray}
&   & (c_2(X)+(n_1+n_2N)F).H \nonumber \\
& = & (f^*c_2(Y)+C_1+n_2NF+(n_1+x)F).H \nonumber \\
& \geq & c_2(Y).f_*H+n_2NF.H \nonumber \\
& = & \lim_{k\to\infty}(a_kC+R_k).f_*H+n_2NF.H \nonumber \\
& = & \lim_{k\to\infty}(a_kf^*C.H+R_k.f_*H+n_2NF.H) \nonumber \\
& = & \lim_{k\to\infty}(a_kC_1.H+(a_ky+n_2N)F.H+R_k.f_*H)\geq 0. \nonumber
\end{eqnarray}
In other words, there exists a positive integer $n$ such that $c_2(X)+nl$ is pseudo-effective.

Let $K=\{u\geq 0|c_2(X)+ul$ is pseudo-effective\}. If $\inf K=0$, then $c_2(X)$ is pseudo-effective. Otherwise, let $r=\inf K>0$. Assume that $\R_+[l]$ is the extremal ray with respect to $f:X\rightarrow Y$. It follows from Lemma \ref{C:2.3} that we can take an ample Cartier divisor $L$ on $X$, a sufficiently small $\varepsilon>0$ and the cone decomposition
\begin{eqnarray}
\overline{NE}(X) & = & \R_+[l]+\sum_{i\in I}\R_+[l_i]+\overline{NE_\varepsilon}(X) \nonumber \\
c_2(X) & = & -rl+\sum_{i\in I}b_il_i+z \label{C:eqn1}
\end{eqnarray}
where $\{\R_+[l_i]\}_{i\in I}$ are extremal rays other than $\R_+[l]$ and $b_i\geq 0$, such that $z.(-K_X)<r/2$. 

If there exists some $\R_+[l_i]$ such that the corresponding extremal contraction $f_i: X\rightarrow Y_i$ is of type ${\rm (F_I)}$ or ${\rm (F_{II})}$, then $c_2(X)$ is pseudo-effective. Otherwise for each $i\in I$, there exists a positive integer $n_i$ such that $c_2(X)+n_il_i$ is pseudo-effective. Note that $c_1(X).c_2(X)=24\chi(\OO_X)=0$. Applying $-K_X$ to each side of (\ref{C:eqn1}), we have
\[ 0=-rl.(-K_X)+\sum_{i\in I}b_il_i.(-K_X)+z.(-K_X), \]
which implies that $\sum b_i>r/12$ since $l.(-K_X)\geq 1$ and $l_i.(-K_X)\leq 6$ (cf.\ \cite{ka91}, Theorem 1).

Let $M=\max_{i\in I}\{n_i\}$. Consider the following pseudo-effective 1-cycle
\begin{eqnarray}
 &      & \sum_{i\in I}b_i(c_2(X)+n_il_i) \nonumber \\
 & \leq & (\sum_{i\in I}b_i)c_2(X)+M(\sum_{i\in I}b_il_i+z) \nonumber \\
 &  =   & (\sum_{i\in I}b_i+M)c_2(X)+Mrl \nonumber
\end{eqnarray}
Since $\sum_{i\in I}b_i>r/12>0$, we have $c_2(X)+Mr/(M+r/12)l$ is pseudo-effective. It contradicts the definition of $r$. Hence $c_2(X)$ is pseudo-effective.

In general, we may take a $\Q$-factorialization $\mu: X^\Q\rightarrow X$ of $X$, namely $X^\Q$ is $\Q$-factorial terminal, and $\mu$ is a projective birational morphism isomorphic in codimension one (cf.\ \cite{ka88}, Corollary 4.5). It is easy to see that $-K_{X^\Q}$ is nef, $\nu(-K_{X^\Q})=2$ and $q(X^\Q)=1$, then $c_2(X^\Q)$ is pseudo-effective, hence so is $c_2(X)=\mu_*c_2(X^\Q)$.
\end{proof}

\begin{rem}\label{C:2.5}
If $q(X)=0$, then $X$ is not necessarily Gorenstein. Provided that every extremal contraction from $X$ is good and $c_1(X).c_2(X)\geq 0$ holds, then the above method is also valid for the case $q(X)=0$. Indeed, let $r_X$ be the Gorenstein index of $X$, then $l.(-K_X)\geq 1/r_X$. If we take $\varepsilon$ to be sufficiently small such that $z.(-K_X)<r/(2r_X)$, then similar arguments work.
\end{rem}

\section{Reduction of the case $q(X)=0$}\label{C:S3}

In order to apply the MMP in the case $q(X)=0$, we should introduce the notion of almost nef divisors, which has been first defined in \cite{ps}.

\begin{defn}\label{C:3.1}
Let $X$ be a normal variety. A $\Q$-Cartier divisor $D$ on $X$ is said to be almost nef, if $D.C\geq 0$ for any irreducible curve $C$ on $X$, except a finite number of rational curves. Such rational curves are said to be exceptional for $D$.
\end{defn}

The following proposition describes the structure of extremal contractions from such regular threefolds whose anticanonical divisors are almost nef.

\begin{prop}\label{C:3.2}
Let $X$ be a $\Q$-factorial terminal projective threefold such that $-K_X$ is almost nef, $\kappa(X)=-\infty$ and $q(X)=0$. Let $f: X\rightarrow Y$ be an extremal contraction. Then $(Y,f)$ is one of the following cases.

${\rm (F_0)}$ $X$ is a Fano threefold with $\rho(X)=1$.

${\rm (F_I)}$ $f$ is a del Pezzo fibration, and $Y\cong\PP^1$.

${\rm (F_{II})}$ $f$ is a conic bundle with discriminant locus $\Delta=\emptyset$, and $Y$ is a canonical surface with $K_Y\equiv 0$.

${\rm (F_{III})}$ $f$ is a conic bundle with discriminant locus $\Delta$, such that $-(4K_Y+\Delta)$ is almost nef, and $Y$ is a rational log terminal surface.

{\rm (D)} $f$ is a divisorial contraction, and $-K_Y$ is almost nef.

{\rm (S)} $f$ is a small contraction. Let $\varphi: X\dashrightarrow X^+$ be the flip of $f$. Then $-K_{X^+}$ is almost nef.
\end{prop}

\begin{proof}
(\ref{C:3.2}.1) If $\dim Y\leq 1$, these are ${\rm (F_0)}$ and ${\rm (F_I)}$.

(\ref{C:3.2}.2) Assume that $\dim Y=2$. Then $Y$ has only quotient singularities, hence is log terminal. Let $S={\rm Sing}X$, $S'=f(S)\subset Y$, $Y_0=Y\setminus S'$, and $X_0=f^{-1}(Y_0)$. Then $f_0: X_0\rightarrow Y_0$ is a usual conic bundle since $X_0$ is smooth. Let $\Delta_0\subset Y_0$ be the discriminant locus of $f_0$ and $\Delta=\overline{\Delta}_0 \subset Y$. It follows from Lemma 1.6 of \cite{ps} that $-(4K_Y+\Delta)$ is almost nef. We divide into two subcases by the emptyness of $\Delta$.

(\ref{C:3.2}.2.1) $\Delta\neq\emptyset$

We claim that $\kappa(Y)=-\infty$. Otherwise, there exists an ample Cartier divisor $H$ on $Y$ such that $K_Y.H\geq 0$. Since $-(4K_Y+\Delta)$ is almost nef, we have $(4K_Y+\Delta).H\leq 0$, hence $0<\Delta.H\leq -4K_Y.H\leq 0$. This is absurd. Thus $p_2(Y)=q(Y)=0$ implies that $Y$ is rational. This is case ${\rm (F_{III})}$.

(\ref{C:3.2}.2.2) $\Delta=\emptyset$

Then $-K_Y$ is almost nef. Let $g: Z\rightarrow Y$ be the minimal resolution of $Y$. We may write $K_Z=g^*K_Y+\sum a_iE_i$, where $E_i$ are exceptional curves, $-1<a_i\leq 0$. Let $C_0$ be any curve on $Z$ such that $C_0\neq E_i$ and $g_*(C_0)$ is not exceptional for $-K_Y$. Then $(-K_Z).C_0=(-K_Y).g_*(C_0)+\sum(-a_i)E_i.C_0\geq 0$. Thus $-K_Z$ is almost nef.

Let $h_1: Z\rightarrow W_1$ be a contraction of a $(-1)$-curve $F_1$. Then $K_Z=h_1^*K_{W_1}+F_1$. Let $C_1$ be any curve on $W_1$ such that $h_{1*}^{-1}(C_1)$ is not exceptional for $-K_Z$, then $(-K_{W_1}).C_1=(-K_Z).h_{1*}^{-1}(C_1)+F_1.h_{1*}^{-1}(C_1)\geq 0$, namely, $-K_{W_1}$ is almost nef. If $\kappa(Y)=\kappa(Z)=-\infty$, then this is case ${\rm (F_{III})}$. Otherwise, there is a birational morphism $h: Z\rightarrow W$, which contracts all $(-1)$-curves such that $K_W$ is nef. Let $L$ be an ample Cartier divisor on $W$, then $-K_W.L\geq 0$ since $-K_W$ is almost nef. But $K_W$ is nef, we only have $K_W\equiv 0$. Therefore $K_Z=h^*K_W+\sum_{j\in J}F_j\equiv\sum_{j\in J}F_j$. By a similar argument, since $-K_Z$ is almost nef, we have that $J=\emptyset$ and $K_Z\equiv 0$.

Let $C$ be any curve on $Y$. Then $K_Y.C=g^*K_Y.g_*^{-1}C=\sum(-a_i)E_i.g_*^{-1}C\geq 0$, namely $K_Y$ is nef. Since $-K_Y$ is almost nef, then $K_Y\equiv 0$. Furthermore, we have $a_i=0$ for all $i$, which implies that $Y$ is canonical. This is case ${\rm (F_{II})}$.

(\ref{C:3.2}.3) If $\dim Y=3$, then $f$ is birational. When $f$ is divisorial, $-K_Y$ is almost nef by Proposition 2.1 of \cite{ps}. When $f$ is small, $-K_{X^+}$ is almost nef by Proposition 2.2 of \cite{ps}.
\end{proof}

For convenience, we give some definitions.

\begin{defn}\label{C:3.3}
Let $\II$=$\bigl\{\Q$-factorial terminal projective threefold $X_0$ $|$ $-K_{X_0}$ is nef, $\nu(-K_{X_0})=2$ and $q(X_0)=0\bigr\}$.

A threefold $X$ is in the class $\A$, if there is a composition of birational maps between some $X_0\in\II$ and $X$: $X_0\dashrightarrow X_1\dashrightarrow\cdots\dashrightarrow X_r=X$, where $\dashrightarrow$ is either a divisorial contraction or a flip.
\end{defn}

Note that for any $X\in\A$, $X$ is a $\Q$-factorial terminal projective threefold such that $-K_X$ is almost nef, $\kappa(X)=-\infty$ and $h^i(\OO_X)=0$ for $i\geq 1$.

We fix the following notation until otherwise stated. Let $f: X\rightarrow Y$ be the extremal contraction induced by an extremal ray $\R_+[l]$ of $X$. The following lemmas show that such extremal contractions are good in some sense.

\begin{lem}\label{C:3.4}
Assume that we are in case ${\rm (F_{II})}$ or ${\rm (F_{III})}$. Then there exists a positive integer $n$ such that $c_2(X)+nl$ is pseudo-effective.
\end{lem}

\begin{proof}
With the same notation as in the proof of Proposition \ref{C:3.2}. Note that $f_0: X_0\rightarrow Y_0$ is a usual conic bundle between smooth quasi-projective varieties, $X\setminus X_0$ is numerically equivalent to some multiple of $l$, and $Y\setminus Y_0$ is a finite number of points. We have the following exact sequence:
\[ 0\rightarrow f_0^*\Omega_{Y_0}\rightarrow \Omega_{X_0}\rightarrow \OO_{X_0}(K_{X_0/Y_0})\rightarrow \OO_{\Gamma_0}\rightarrow 0 \]
where $\Gamma_0$ is defined as in Lemma 3.6 of \cite{xie}. Then
\begin{eqnarray}
c_2(X_0) & = & f_0^*c_2(\Omega_{Y_0})+f_0^*c_1(\Omega_{Y_0}).K_{X_0/Y_0}+\Gamma_0 \nonumber \\
         & = & f_0^*(c_2(Y_0)-c_1^2(Y_0))+f_0^*c_1(Y_0).c_1(X_0)+\Gamma_0. \nonumber
\end{eqnarray}
Therefore $c_2(X)=f^*(c_2(Y)-c_1^2(Y))+f^*c_1(Y).c_1(X)+\Gamma+al$, where $\Gamma=\overline{\Gamma}_0$, and $a$ is a rational number.

${\rm (F_{II})}$ There exists a positive integer $n$ such that $c_2(X)+nl$ is pseudo-effective, since $K_Y\equiv 0$.

${\rm (F_{III})}$ We rewrite the above formula as follows:
\begin{eqnarray}
c_2(X)+nl=n_0l+(1/4)f^*(-4K_Y-\Delta).(-K_X)+(1/4)f^*\Delta.(-K_X)+\Gamma \label{C:eqn2}
\end{eqnarray}
where $n\in\N$, and $n_0\in\Q^+$. For an ample Cartier divisor $H$ on $X$, assume that $mH$ is very ample for some $m\in\N$. Since $-(4K_Y+\Delta)$ and $-K_X$ are almost nef, we may choose a general member $M\in |mH|$ such that $M$ contains no exceptional curves for $-K_X$ and contains no curves in $\cup_\lambda\Supp f^*C_\lambda.(-K_X)$, where \{$C_\lambda$\} is the set of exceptional curves for $-(4K_Y+\Delta)$. It is easy to show that $c_2(X)+nl$ is pseudo-effective by applying $M$ to each side of (\ref{C:eqn2}).
\end{proof}

\begin{lem}\label{C:3.5}
Assume that we are in case {\rm (D)}. Furthermore, assume that $c_2(Y)$ is pseudo-effective. Then there exists a positive integer $n$ such that $c_2(X)+nl$ is pseudo-effective.
\end{lem}

\begin{proof}
(\ref{C:3.5}.1) $f$ contracts a divisor $E$ to a point $p\in Y$.

Because $f: X\setminus E\rightarrow Y\setminus \{p\}$ is an isomorphism, and both $X$ and $Y$ are of terminal singularities, we have $c_2(X\setminus E)=f^*c_2(Y\setminus \{p\})=f^*c_2(Y)$. Thus $c_2(X)$ and $f^*c_2(Y)$ differ by a 1-cycle whose support is contained in $E$. We may write
\[ c_2(X)=f^*c_2(Y)+\sum\alpha_il_i, \quad \alpha_i\in\Q, l_i\subset E. \]
Since $f_*(l_i)\subset f_*(E)=p$, some positive multiple of each $l_i$ is numerically equivalent to $l$. Since $c_2(Y)$ is pseudo-effective, there exists a positive integer $n$ such that $c_2(X)+nl$ is pseudo-effective.

(\ref{C:3.5}.2) $f$ contracts a divisor $E$ to a curve $C$ in $Y$.

Let $S={\rm Sing}X$, $S'=f(S)\subset Y$, $Y_0=Y\setminus S'$, and $X_0=f^{-1}(Y_0)$. Then $f_0: X_0\rightarrow Y_0$ is a blow-up along a smooth curve $C$. By a similar argument to that in Theorem \ref{C:2.4}, we can prove that there exists a positive integer $n$ such that $c_2(X)+nl$ is pseudo-effective.
\end{proof}

\begin{lem}\label{C:3.6}
Assume that we are in case {\rm (S)}. Then $c_2(Y)$ is pseudo-effective if and only if $c_2(X)+nl$ is pseudo-effective for some $n\in\N$.
\end{lem}

\begin{proof}
Since $f$ is small, the exceptional locus $E$ consists of finitely many rational curves, and $f(E)$ is a finite set of points. Since $f: X\setminus E\rightarrow Y\setminus f(E)$ is an isomorphism, we have $f_*c_2(X\setminus E)=c_2(Y\setminus f(E))=c_2(Y)$. On the other hand, $c_2(X\setminus E)=c_2(X)+\alpha l$ for some $\alpha\in\Q$. Then we have $f_*c_2(X)=c_2(Y)$.

For the ``if'' part: $f_*(c_2(X)+nl)=f_*c_2(X)=c_2(Y)$ is pseudo-effective.

For the ``only if'' part: $c_2(Y)$ is pseudo-effective, so is $c_2(X\setminus E)$, hence so is $c_2(X)+nl$ for some $n\in\N$.
\end{proof}

\begin{lem}\label{C:3.7}
Assume that we are in case ${\rm (F_I)}$. Then $c_2(X)$ is pseudo-effective provided that $c_1(X).c_2(X)\geq 0$ holds.
\end{lem}

\begin{proof}
Let $X_\xi=f^{-1}(\xi)$ for a general point $\xi\in Y=\PP^1$. Since $\Pic(X)\cong f^*\Pic(\PP^1)\oplus\Z$, for any ample Cartier divisor $M$ on $X$, we may write $M\equiv a(-K_X)+bX_\xi$ for some $a,b\in\Q$. Hence $M.l=a(-K_X).l>0$, which shows that $a>0$.

(\ref{C:3.7}.1) If $-K_X$ is nef and $\nu(-K_X)=1$ or 3, then $c_2(X)$ is pseudo-effective by Theorem 2.2 of \cite{xie}.

(\ref{C:3.7}.2) If $-K_X$ is nef and $\nu(-K_X)=2$, then $(-K_X)^3=0$. Hence $0<M.(-K_X)^2=bX_\xi.(-K_X)^2=b(-K_{X_\xi})^2$, which shows that $b>0$.

(\ref{C:3.7}.3) If there exists an exceptional curve $C$ for $-K_X$, then $f(C)$ is not a point since $f$ is an extremal contraction. Therefore $f: C\rightarrow \PP^1$ is surjective, hence $X_\xi.C>0$. We have $bX_\xi.C=M.C+aK_X.C>0$, which shows that $b>0$.

Thus it is sufficient to prove that $c_1(X).c_2(X)\geq 0$ and $X_\xi.c_2(X)\geq 0$ for verifying the pseudo-effectivity of $c_2(X)$. $X_\xi.c_2(X)=c_2(X_\xi)>0$ since $X_\xi$ is a smooth del Pezzo surface. In the subcase (\ref{C:3.7}.2), since $-K_X$ is nef, we always have $c_1(X).c_2(X)\geq 0$ (cf.\ \cite{xie}, Theorem 2.2). In the subcase (\ref{C:3.7}.3), this follows from the assumption.
\end{proof}

It is natural to put forward the following problem, and we will make some discussions on this problem in the end of this section.

\begin{prob}\label{C:3.8}
Does $c_1(X).c_2(X)\geq 0$ hold for any $X\in\A$?
\end{prob}

We may reduce the case $q(X)=0$ to Problem \ref{C:3.8} by the following:

\begin{thm}\label{C:3.9}
$c_2(X)$ is pseudo-effective for any $X\in\A$ provided that Problem \ref{C:3.8} is true. In particular, let $X_0$ be a terminal projective threefold such that $-K_{X_0}$ is nef, $\nu(-K_{X_0})=2$ and $q(X_0)=0$. Then $c_2(X_0)$ is pseudo-effective provided that $c_1(Y).c_2(Y)\geq 0$ for any $Y\in\A(X^\Q_0)$, where $\mu: X^\Q_0\rightarrow X_0$ is a $\Q$-factorialization of $X_0$.
\end{thm}

\begin{proof}
We use induction on the Picard number $\rho(X)$.

It is easy to see that $c_2(X)$ is pseudo-effective when $\rho(X)=1$, since only ${\rm (F_0)}$ occurs. Assume that the conclusion holds for $\rho(X)<\rho$. Let $\rho(X)=\rho\geq 2$.

We recall the definition of difficulty (cf.\ \cite{sh}). Let $g: \widetilde{X}\rightarrow X$ be a resolution of $X$. We may write $K_{\widetilde{X}}=g^*K_X+\sum a_iE_i$, where $E_i$ are the exceptional divisors of $g$. The difficulty of $X$ is defined to be $d(X)=\sharp\{i$ $|$ $a_i<1\}$. Note that $d(X)$ is independent of the choice of the resolution $g$, hence is well-defined. In fact, $d(X)$ has been introduced to prove the termination of the flips.

We use induction on $d(X)$. If $d(X)=0$, then case (S) cannot occur. Otherwise, the existence of the flip $\varphi: X\dashrightarrow X^+$ yields $d(X^+)<d(X)=0$, this is absurd. In the other cases, $c_2(X)+nl$ is pseudo-effective for some $n\in\N$ by the preceding lemmas, and $c_1(X).c_2(X)\geq 0$ by Problem \ref{C:3.8}. Hence by a similar argument to that of Theorem \ref{C:2.4} (cf.\ Remark \ref{C:2.5}), we can prove that $c_2(X)$ is pseudo-effective.

Assume that $c_2(X)$ is pseudo-effective for $d(X)<d$. Let $d(X)=d$. In case (S), note that $X^+\in\A$, $\rho(X^+)=\rho(X)=\rho$, and $d(X^+)<d(X)=d$, so $c_2(X^+)$ is pseudo-effective by induction hypothesis. Therefore $c_2(Y)=f^+_*c_2(X^+)$ is pseudo-effective since $f^+: X^+\rightarrow Y$ contracts several rational curves to points. It follows from Lemma \ref{C:3.6} that $c_2(X)+nl$ is pseudo-effective for some $n\in\N$. By the same argument, we can prove that $c_2(X)$ is pseudo-effective.
\end{proof}

As a special case, we give the following:

\begin{cor}\label{C:3.10}
Let $X\in\A$. If $X$ is factorial and $\rho(X)\leq 3$, then $c_2(X)$ is pseudo-effective.
\end{cor}

\begin{proof}
(\ref{C:3.10}.1) If $\rho(X)=1$, then the pseudo-effectivity of $c_2(X)$ is obvious.

(\ref{C:3.10}.2) If $\rho(X)=2$, then we can prove a stronger statement.

{\bf Claim.} {\it Let $X\in\A$. Assume that $\rho(X)=2$ and that $X$ has at most one cyclic quotient singularity of type $\frac{1}{2}(1,1,1)$ and the other singular points are Gorenstein. Then $c_2(X)$ is pseudo-effective.}

{\it Proof of the Claim.} We use induction on the difficulty $d(X)$. If $d(X)=0$, then case (S) cannot occur. In the other cases, it is easy to show that $c_2(X)$ or $c_2(X)+nl$ is pseudo-effective for some $n\in\N$. Since $c_1(X).c_2(X)>0$, we can prove that $c_2(X)$ is pseudo-effective. Assume that $c_2(X)$ is pseudo-effective for $d(X)<d$. Let $d(X)=d$. It is sufficient to show that in case (S), $X^+$ also satisfies the assumption of the Claim. Indeed, let $q\in X^+$ be a terminal singular point of index $r>1$. If $q$ is not contained in the exceptional locus of $f^+: X^+\rightarrow Y$, then $q\in X^+$ is a cyclic quotient singularity of type $\frac{1}{2}(1,1,1)$. Otherwise, let $\widetilde X$ be a common resolution of $X$ and $X^+$. Then there is an exceptional divisor $E$ on $\widetilde X$ such that $a(E,X^+)=1/r$ (cf.\ \cite{ka92}). By Lemma 3.38 of \cite{km}, we have $1/2\leq a(E,X)<a(E,X^+)=1/r$, hence $r<2$. The contradiction completes the proof of the Claim.

(\ref{C:3.10}.3) If $\rho(X)=3$, then case (S) cannot occur since $X$ is factorial. In case (D), note that $Y$ satisfies the conditions of the Claim (cf.\ \cite{cu}), so $c_2(Y)$ is pseudo-effective. Since $c_1(X).c_2(X)>0$, we can prove that $c_2(X)$ is pseudo-effective.
\end{proof}

It follows from the singular Riemann-Roch formula (cf.\ \cite{re87}, Corollary 10.3) that
\[ 1=\chi(\OO_X)=\frac {1}{24}c_1(X).c_2(X)+\frac {1}{24}\sum_i(r_i-\frac {1}{r_i}), \]
where $\{\frac {1}{r_i}(a_i,-a_i,1), (a_i,r_i)=1 \}_i$ are the fictitious cyclic quotient singularities for $X$.

Define $F(X)=\sum_i(r_i-1/r_i)$ to be the fictitious number of $X$. Let $f: X\rightarrow Y$ be a divisorial contraction. If $f$ contracts a divisor to a point (resp.\ a curve), we say that $f$ is of type O (resp.\ I). If $\varphi: X\dashrightarrow X^+$ is a flip, we say that $\varphi$ is of type F.

In general, when we run the MMP, the latter variety should have better singularities than the former one. For a terminal threefold $X$, $F(X)$ may reflect some information of its singularities. We expect that $F(X)$ decreases in type I and F. In type O, there is an example to show that $F(X)$ may increase strictly.

\begin{ex}\label{C:3.11}
Let $p\in Y$ be a germ of cyclic quotient singularity of type $\frac{1}{r}(a,-a,1)$, where $(r,a)=1,r>1$. Let $f: X\rightarrow Y$ be the weighted blow-up along $p$ with weight $\frac{1}{r}(a,r-a,1)$. Then $X$ has two cyclic quotient singular points of type $\frac{1}{a}(r,-r,1)$ and $\frac{1}{r-a}(a,-a,1)$ respectively. It is easy to see that $F(X)<F(Y)$.
\end{ex}

In order to work out Problem \ref{C:3.8}, it is helpful to consider the following:

\begin{prob}\label{C:3.12}
Let $X$ be a $\Q$-factorial terminal projective threefold, and $\varphi: X\dashrightarrow X^+$ a flip. Does $F(X)\geq F(X^+)$ hold?
\end{prob}

\begin{prob}\label{C:3.13}
Let $X$ be a $\Q$-factorial terminal projective threefold, and $f: X\rightarrow Y$ a divisorial contraction of type I. Does $F(X)\geq F(Y)$ hold?
\end{prob}

We give an example as an evidence of Problem \ref{C:3.12}.

\begin{ex}
With the same notation as in Problem \ref{C:3.12}. If $\varphi: X\dashrightarrow X^+$ has a deformation consisting of toric flips, then $F(X)>F(X^+)$ holds. Indeed, since $F(X)$ is additive under the deformations, we may assume that $\varphi: X\dashrightarrow X^+$ is a toric flip. By the structure theorem for toric morphisms (cf. \cite{re83a}), a toric flip is given by the following formula and diagram:
\[ a_1e_1+a_2e_2=a_3e_3+a_4e_4, \]
\vskip 5mm
\unitlength=1mm
\begin{center}
$X:$\qquad
\begin{picture}(22,10)(0,0)
\put(-4,0){$e_3$}
\put(0,0){\line(1,1){10}}
\put(0,0){\line(1,-1){10}}
\put(10,10){\line(1,-1){10}}
\put(10,10){\line(0,-1){20}}
\put(9,11){$e_1$}
\put(10,-10){\line(1,1){10}}
\put(9,-13){$e_2$}
\put(21,0){$e_4$}
\end{picture}
\quad
$\dashrightarrow$
$X^+:$\qquad
\begin{picture}(22,10)(0,0)
\put(-4,0){$e_3$}
\put(0,0){\line(1,1){10}}
\put(0,0){\line(1,0){20}}
\put(0,0){\line(1,-1){10}}
\put(10,10){\line(1,-1){10}}
\put(9,11){$e_1$}
\put(10,-10){\line(1,1){10}}
\put(9,-13){$e_2$}
\put(21,0){$e_4$}
\end{picture}
\end{center}
\vskip 10mm

It follows from the White-Frumkin theorem (cf. \cite{da}) that there are only two types of toric flips in dimension three.

Let \{$e_1,e_2,e_3$\} be a basis of the lattice $N\cong\Z^3$.

(i) $ae_1+(r-a)e_2=re_3+e_4$, where $(r,a)=1, r>a$. $X$ has one cyclic quotient singular point of type $\frac{1}{r}(a,-a,1)$, $X^+$ has two cyclic quotient singular points of type $\frac{1}{r-a}(a,1,-a)$ and $\frac{1}{a}(1,r,-r)$ respectively.

(ii) $ae_1+e_2=re_3+e_4$, where $(r,a)=1, r>a$. $X$ has one cyclic quotient singular point of type $\frac{1}{r}(a,-a,1)$, $X^+$ has one cyclic quotient singular point of type $\frac{1}{a}(1,r,-r)$.

It is easy to see that $F(X)>F(X^+)$ holds for each case.
\end{ex}

\section{Applications}\label{C:S4}

The following Effective Non-vanishing Conjecture has been put forward by Ambro and Kawamata (cf.\ \cite{am,ka00}).

\begin{conj}\label{C:4.1}
Let $X$ be a complete normal variety, $B$ an effective $\R$-divisor on $X$ such  that the pair $(X,B)$ is Kawamata log terminal (KLT, for short), and $D$ a Cartier divisor on $X$. Assume that $D$ is nef and that $D-(K_X+B)$ is nef and big. Then $H^0(X,D)\neq 0$.
\end{conj}

\begin{rem}\label{C:4.2}
(\ref{C:4.2}.1) By the Kawamata-Viehweg vanishing theorem, we have $H^i(X,D)=0$ for any positive integer $i$. Thus $H^0(X,D)\neq 0$ is equivalent to $\chi(X,D)\neq 0$.

(\ref{C:4.2}.2) When $\dim X=1$, the effective non-vanishing follows easily 
from the Riemann-Roch theorem. When $\dim X=2$, Conjecture \ref{C:4.1} 
has been proven by Kawamata by means of the semipositivity theorem 
(cf.\ \cite{ka00}). For higher dimensional cases, Conjecture \ref{C:4.1} 
is still open.

(\ref{C:4.2}.3) If $X$ is a toric variety, then Conjecture \ref{C:4.1} is trivial. Indeed, any nef invariant Cartier divisor is base point free (cf.\ \cite{mu}, Theorem 3.1), hence is non-vanishing.
\end{rem}

In fact, we can prove that $c_2(X)$ being pseudo-effective is a sufficient condition such that the Effective Non-vanishing Conjecture holds for terminal projective threefolds with $B=0$.

\begin{prop}\label{C:4.3}
Let $X$ be a terminal projective threefold such that $c_2(X)$ is pseudo-effective. Let $D$ be a Cartier divisor on $X$ such that $D$ is nef and $D-K_X$ is nef and big. Then $H^0(X,D)\neq 0$.
\end{prop}

\begin{proof}
If $\nu(D)\leq 2$, then by Theorem 2.2 of \cite{ka00}, we may reduce this case to the log surface case. So we may assume that $D$ is nef and big. By the singular Riemann-Roch formula, we have
\begin{eqnarray}
&   & h^0(\OO_X(D)) \nonumber \\
& = & \frac {1} {12}D(D-K_X)(2D-K_X)+\frac {1} {12}D.c_2(X)+\chi(\OO_X) \nonumber \\
& \geq & \frac {1} {12}D^2(D-K_X)+\frac {1} {12}D(D-K_X)^2+\frac {1} {24}(2D-K_X).c_2(X)>0. \nonumber
\end{eqnarray}
\end{proof}

We have the following two applications, where Corollary \ref{C:4.4} was first proved by Kawamata (cf.\ \cite{ka00}, Proposition 4.1).

\begin{cor}\label{C:4.4}
Let $X$ be a canonical projective minimal threefold. 
Then the effective non-vanishing holds on $X$.
\end{cor}

\begin{cor}\label{C:4.5}
Let $X$ be a canonical projective threefold such that $-K_X$ is nef and either

{\rm (i)} $\nu(-K_X)\neq 2$, or

{\rm (ii)} $\nu(-K_X)=2$ and $q(X)=1$.

\noindent Then the effective non-vanishing holds on $X$.
\end{cor}

\begin{proof}[Proof of Corollaries \ref{C:4.4} and \ref{C:4.5}]
It follows from the main theorem of \cite{re83b} that there exists a partial resolution $g: X'\rightarrow X$ such that $X'$ is terminal and $K_{X'}=g^*K_X$. 
Thus we may replace $X$ by $X'$ for proving such non-vanishings. 
For Corollary \ref{C:4.4}, the rest is due to the Miyaoka theorem 
(cf.\ \cite{mi}). 
For Corollary \ref{C:4.5}, the rest is due to Theorems \ref{C:1.2} 
and \ref{C:2.4}.
\end{proof}

\textsc{Graduate School of Mathematical Sciences, University of Tokyo, 3-8-1 Komaba, Meguro, Tokyo 153-8914, Japan}

\textit{E-mail address}: \texttt{xqh@ms.u-tokyo.ac.jp}

\end{document}